\documentclass{amsart}
\usepackage[final]{epsfig}
\usepackage{graphics}
\usepackage{float}
\usepackage{amsmath}
\usepackage{amsfonts}
\usepackage{latexsym}
\usepackage{amssymb}
\usepackage{graphicx}
\usepackage{url}
\usepackage{epstopdf}
\usepackage{verbatim}
\usepackage[margin=1in]{geometry}
\usepackage{amsthm}
\usepackage{tikz}
\usepackage{caption}
\usepackage{enumerate}

\makeatletter
\newtheorem*{rep@theorem}{\rep@title}

\newcommand{\newreptheorem}[2]{%
\newenvironment{rep#1}[1]{%
 \def\rep@title{#2 \ref{##1}}%
 \begin{rep@theorem}}%
 {\end{rep@theorem}}}
\makeatother

\newcommand{\thmrestatement}[1] {\vspace{2 mm} \noindent \textbf{Theorem #1. }\vspace{2 mm}}

\newtheorem{lemma}{Lemma}[section]

\newtheorem{example}[lemma]{Example}
\newtheorem{theorem}[lemma]{Theorem}

\newtheorem*{theorem*}{Theorem}

\newreptheorem{theorem}{Theorem}

\newcommand{\proofend}{$\Box$\bigskip}

\newcommand{\N}{{\mathbb N}}

\newcommand{\Z}{{\mathbb Z}}

\newcommand{\lfl}{\lfloor}
\newcommand{\rfl}{\rfloor}
\newcommand{\lrep}{\left[}
\newcommand{\rrep}{\right]}
\def\proof{\par{\it Proof}. \ignorespaces}

\DeclareMathOperator{\Aut}{Aut}

\begin{document}

\title{Equality of Dedekind sums mod $8 \mathbb{Z}$}

\author{Emmanuel Tsukerman}
\date{\today}
\address{Department of Mathematics, University of California,
Berkeley, CA 94720-3840}
\email{e.tsukerman@berkeley.edu}
\subjclass[2000]{}

\begin{abstract}
Using a generalization due to Lerch [M. Lerch, Sur un th\'{e}or\`{e}me de Zolotarev. Bull. Intern. de l'Acad. Fran\c{c}ois Joseph 3 (1896), 34-37] of a classical lemma of Zolotarev, employed in Zolotarev's proof of the law of quadratic reciprocity, we determine necessary and sufficient conditions for the difference of two Dedekind sums to be in $8\mathbb{Z}$. These yield new necessary conditions for equality of two Dedekind sums. In addition, we resolve a conjecture of Girstmair [Girstmair, 
Congruences mod 4 for the alternating sum of the partial quotients, arXiv: 1501.00655].
\end{abstract}

\maketitle


{\bf \small MSC:} 11F20.
\[\]

\section{Background}

Dedekind sums  are classical objects of study  introduced by Richard Dedekind in the 19th century in his study of the $\eta$-function \cite{dedekind}. Among many other areas of mathematics, Dedekind sums appear in: geometry (lattice point enumeration in polytopes \cite{MR2271992}), topology (signature defects of manifolds \cite{MR0650832}) and algorithmic complexity (pseudo random number generators \cite{Knuth:1997:ACP:270146}). To define the Dedekind sums, let

\[
((x))=\begin{cases}
x-\lfloor x\rfloor - 1/2, &\mbox{if }x\in\mathbb{R}\setminus\mathbb{Z};\\
0,&\mbox{if }x\in\mathbb{Z}.
\end{cases}
\]
Then the Dedekind sum $s(a,b)$ for $a,b \in \N$ coprime is defined by
\[
s(a,b)=\sum_{k=1}^b ((\frac{ak}{b}))((\frac{k}{b})).
\]

 Recently, in \cite{MR2873148}, Jabuka et al. raise the question of when two Dedekind sums $s(a_1,b)$ and $s(a_2,b)$ are equal. In the same paper, they prove the necessary condition $b \mid (a_1 a_2-1)(a_1-a_2)$ for equality of two dedekind sums $s(a_1,b)$ and $s(a_2,b)$. Girstmair \cite{MR3189994} shows that this condition is equivalent to $12s(a_1,b)-12s(a_2,b) \in \mathbb{Z}$. In \cite{doi:10.1142/S1793042115500748}, necessary and sufficient conditions for $12s(a_1,b)-12s(a_2,b) \in  2\mathbb{Z}, 4\mathbb{Z}$ are given. 

In this note we give necessary and sufficient conditions for $12s(a_1,b)-12s(a_2,b) \in 8\mathbb{Z}$ by using a generalization of Zolotarev's classical lemma relating the Jacobi symbol to the sign of a special permutation\footnote{The motivation behind Zolotarev's work was to produce a proof of the law of quadratic reciprocity.} due to Lerch \cite{Lerch}. Along the way, we resolve a conjecture of Girstmair about the alternating sum of partial quotients modulo 4 \cite{GirstAlt}.

\section{Preliminaries}

Let  $ \pi_{(a,b)} \in \Aut(\Z / b \Z), \pi_{(a,b)}:x \mapsto ax$. Let $\lrep x \rrep_b=x-b \lfl \frac{x}{b} \rfl$ be the function taking $x \in \Z/ b \Z$ to its smallest nonnegative representative. We view $\pi_{(a,b)}$ as a permutation of $\{0,1,\ldots,b-1\}$ given by
\[
\pi_{(a,b)}=\left(\begin{array}{ccccccc}
0 & 1 & \cdots & b-1 \\
\lrep \pi(0) \rrep & \lrep \pi(1) \rrep & \cdots & \lrep \pi(b-1) \rrep
\end{array}\right)=\left(\begin{array}{ccccccc}
0 & 1 & \cdots & b-1 \\
0& \lrep a \rrep_b & \cdots & \lrep (b-1)a \rrep_b
\end{array}\right).
\]
The precedent for doing so is already present in the work of Zolotarev, in which he relates the sign of $\pi_{(a,b)}$ to the Jacobi symbol and obtains a proof of the law of quadratic reciprocity (see, e.g., \cite[pg. 38]{MR0357299}). Let $I(a,b)$ denote the number of inversions  of $\pi_{(a,b)}$.

\begin{theorem}  \label{zol}(Zolotarev) For odd $b$ and $(a,b)=1$,
\[
(-1)^{I(a,b)}=\left( \frac{a}{b} \right).
\]
\end{theorem}

The following result shows that the inversions of $\pi_{(a,b)}$ and Dedekind sums are closely related.

\begin{theorem} \label{invded}(Meyer, \cite{Meyer1957})  The number of inversions $I(a,b)$ of $\pi_{(a,b)}$ is equal to
\[
I(a,b)=-3b s(a,b)+\frac{1}{4}(b-1)(b-2),
\]
where $s(a,b)$ is the Dedekind sum.
\end{theorem}

From the reciprocity law of dedekind sums, one obtains a reciprocity law for inversions. 

\begin{theorem} \label{sal}(Sali\'{e}, \cite[p. 163]{Meyer1957})  For all coprime $a,b \in \N$
\begin{align}\label{recip}
4aI(a,b)+4bI(b,a)=(a-1)(b-1)(a+b-1).
\end{align}
\end{theorem}

Let $a$ and $b$ be positive integers, $a<b$. Consider the regular continued fraction expansion
\[
\frac{a}{b}=[0,a_1,\ldots,a_n],
\]
where all digits $a_1,\ldots,a_n$ are positive integers. We assume that $n$ is odd\footnote{If $n$ is even, we can consider instead $[0,a_1,\ldots,a_n-1,1]$.}. We will be interested in
\[
T(a,b)=\sum_{j=1}^n (-1)^{j-1} a_j
\]
and
\[
D(a,b)=\sum_{j=1}^n a_j.
\]
With this notation,
\begin{theorem}\label{BHK}(Barkan-Hickerson-Knuth formula) Let $a,b \in \N$ be coprime and let $a^* a \equiv 1 \pmod{b}$ with $0<a^*<b$. Then
\[
12s(a,b)=T(a,b)+\frac{a+a^*}{b}-3.
\]
\end{theorem}

In \cite{Lerch}, Lerch improves upon Zolotarev's Lemma by determining the parity of $I(a,b)$ when $b$ is even:

\begin{theorem}\label{GenZol}(Lerch)
\[
I(a,b) \equiv \begin{cases}
\frac{1-\left( \frac{a}{b} \right)}{2}, & \text{ if }b\text{ is odd} \\
\frac{(a-1)(b+a-1)}{4}, & \text{ if }b\text{ is even}
\end{cases} \pmod{2}.
\]
\end{theorem}

\proof
We assume that $b$ is even, as the result for $b$ odd follows from Theorem \ref{zol}. Reducing the equality
\[
4aI(a,b)+4bI(b,a)=(a-1)(b-1)(a+b-1)
\]
from Theorem \ref{sal} modulo $8$ and using the assumption that $b$ is even yields
\[
4aI(a,b) \equiv (a-1)(b-1)(a+b-1) \pmod{8}.
\]
Since $a-1$ and $a+b-1$ are even, 
\[
aI(a,b) \equiv (b-1) \frac{(a-1)(b+a-1)}{4} \pmod{2},
\]
from which the claim follows.
\proofend

For further generalizations of Zolotarev's Lemma, see \cite{MRZol}.

\section{Main Results}

As a consequence of Theorem \ref{GenZol}, we are able to show the following necessary and sufficient conditions for equality of Dedekind sums mod $8 \mathbb{Z}$.

\begin{theorem} \label{necCond} Let $a_1,a_2 \in \N$ be relatively prime to $b \in \N$. The following are equivalent:
\begin{enumerate}[(a)]
\item \label{a} $I(a_1,b) \equiv I(a_2,b) \pmod{2b}$
\item \label{b} $3s(a_1,b)-3s(a_2,b) \in 2\Z$
\item \label{c} Let $\left( \frac{a}{b} \right)$ denote the Jacobi Symbol and define 
\[
\mu(a,b)=\begin{cases}
\frac{1-\left( \frac{a}{b} \right)}{2}, & \text{ if }b\text{ is odd} \\
\frac{(a-1)(b+a-1)}{4}, & \text{ if }b\text{ is even}
\end{cases} 
\]
Then
\[
(a_1-a_2)(b-1)(b+a_1 a_2-1) \equiv 4b(a_2 \mu(b,a_1)-a_1 \mu(b,a_2)) \pmod{8b}.
\]
\end{enumerate}
\end{theorem}

We also determine $T(a,b)$ mod $8$:

\begin{theorem}\label{altSum} Let $a,b \in \N$ be coprime. Then
\[
b T(a,b) \equiv -4\mu(a,b)+b^2+2-a-a^* \pmod{8}.
\]
\end{theorem}

Reducing further to mod $4$ and mod $2$ resolves a conjecture of Girstmair \cite{GirstAlt}.

\section{Proofs and Examples}

\thmrestatement{\ref{necCond}}
 Let $a_1,a_2 \in \N$ be relatively prime to $b \in \N$. The following are equivalent:
\begin{enumerate}[(a)]
\item  $I(a_1,b) \equiv I(a_2,b) \pmod{2b}$
\item  $3s(a_1,b)-3s(a_2,b) \in 2\Z$
\item  Let $\left( \frac{a}{b} \right)$ denote the Jacobi Symbol and define 
\[
\mu(a,b)=\begin{cases}
\frac{1-\left( \frac{a}{b} \right)}{2}, & \text{ if }b\text{ is odd} \\
\frac{(a-1)(b+a-1)}{4}, & \text{ if }b\text{ is even}
\end{cases} 
\]
Then
\begin{align} \label{cond}
(a_1-a_2)(b-1)(b+a_1 a_2-1) \equiv 4b(a_2 \mu(b,a_1)-a_1 \mu(b,a_2)) \pmod{8b}.
\end{align}
\end{enumerate}

\proof 
The equivalence of \ref{necCond}(\ref{a}) and \ref{necCond}(\ref{b}) follows from Theorem \ref{invded}. Reducing equation (\ref{recip}) of Theorem \ref{sal}  modulo $8b$ and using Theorem \ref{GenZol} yields
\[
4aI(a,b)+4b \mu(b,a) \equiv (a-1)(b-1)(a+b-1) \pmod{8b}.
\]
\proofend

That Theorem \ref{necCond} is not a sufficient condition for the equality of two Dedekind sums is demonstrated in the following example.

\begin{example}
Take $a_1=1, a_2=15$ and $b=49$. Then
\[
\left( \frac{b}{a_1} \right)=1, \quad \left( \frac{b}{a_2} \right)=1.
\]
We have
\[
(a_1-a_2)(b-1)(b+a_1a_2-1)=-42336=108 \cdot 8 \cdot 49 \equiv 0 \pmod{8b}.
\]
Thus we expect $3s(a_1,b)-3s(a_2,b) \in 2 \Z$. Indeed,
\[
s(a_1,b)=\frac{188}{49}, \quad s(a_2,b)=-\frac{8}{49},
\]
so that
\[
3s(a_1,b)-3s(a_2,b)=12.
\]
Equality does not hold.
\end{example}

\thmrestatement{\ref{altSum}} Let $a,b \in \N$ be coprime. Then
\begin{align} \label{Tmod8}
b T(a,b) \equiv -4\mu(a,b)+b^2+2-a-a^* \pmod{8}.
\end{align}

\proof
By Theorems \ref{invded} and \ref{BHK}, we have
\[
bT(a,b)=12bs(a,b)-a-a^*+3b
\]
\[
=-4I(a,b)+b^2+2-a-a^*.
\]
Reducing modulo $8$ and using Theorem \ref{GenZol},
\[
b T(a,b) \equiv -4\mu(a,b)+b^2+2-a-a^* \pmod{8}.
\]

\proofend

Let $k \in \Z$ satisfy $a a^*=1+kb$. In \cite{GirstAlt}, Girstmair conjectures that if $a \equiv a^* \equiv 0 \pmod{2}$, then
\begin{enumerate}[(i)]
\item \label{parta} If $a$ or $a^*$ is $\equiv 2 \pmod4$, then $T(a,b) \equiv (b-k)/2 \pmod4$
\item \label{partb}If $a$ and $a^*$ are both $\equiv 0 \pmod4$, then $T(a,b) \equiv (k-b)/2 \pmod4$
\item \label{partc} If $a$ and $a^*$ are both $\equiv 0 \pmod{4}$, then $D(a,b)$ is odd.
\end{enumerate}
We now show how this follows from Theorem \ref{altSum}. Reducing congruence (\ref{Tmod8}) mod $4$ gives
\[
b T(a,b) \equiv b^2+2-a-a^* \pmod{4}.
\]
Assume first that $a \equiv a^* \equiv 0 \pmod{4}$. Then 
\[
bT(a,b) \equiv b^2+2 \equiv -1 \pmod{4} \implies T(a,b) \equiv -b^{-1} \equiv -b \pmod{4}.
\]
On the other hand, 
\[
1+kb \equiv 0 \pmod{8} \implies k \equiv -b \pmod{8}.
\]
This proves (\ref{partb}). As Girstmair notes, part (\ref{partc}) follows from (\ref{partb}). 

Next we show (\ref{parta}). It suffices to prove the result when $a \equiv 2 \pmod{4}$, since $T(a,b)=T(a^*,b)$. We have
\[
b T(a,b) \equiv 1-a^* \pmod{4} \implies T(a,b) \equiv b^{-1}(1-a^*) \equiv b(1-a^*) \pmod{4}.
\]
On the other hand,
\[
\frac{b-k}{2} \equiv \frac{b-b^{-1}(aa^*-1)}{2} \equiv b-\frac{baa^*}{2} \equiv b-b a (\frac{a^*}{2}) \equiv b-ba^* \pmod{4}.
\]
This completes the proof.

This, together with the results in \cite{GirstAlt}, determines $T(a,b)$ and $D(a,b)$ in all cases.

\bigskip
{\bf Acknowledgments}.  
The author would like to thank Kurt Girstmair and Pete Clark for helpful comments and bibliographic references.

This material is based upon work supported by the National Science Foundation Graduate Research Fellowship under Grant No. DGE 1106400. Any opinion, findings, and conclusions or recommendations expressed in this material are those of the authors(s) and do not necessarily reflect the views of the National Science Foundation.

\bibliographystyle{alpha}
\bibliography{bibliography}

\end{document}